\documentclass[a4paper,11pt]{amsart}

\usepackage {amssymb}
\usepackage{amscd}

\theoremstyle{plain}
\newtheorem{thm}{Theorem}[section]

\newtheorem{cor}[thm]{Corollary}

\def\SL{\operatorname{SL}}
\def\SO{\operatorname{SO}}
\def\Sp{\operatorname{Sp}}
\def\CSp{\operatorname{CSp}}

\def\PSL{\operatorname{PSL}}

\def\Spin{\operatorname{Spin}}

\def\typeA{\mathrm{A}}
\def\typeB{\mathrm{B}}
\def\typeC{\mathrm{C}}

\def\typeF{\mathrm{F}}
\def\typeG{\mathrm{G}}

\newcommand{\FF}{{\mathbb{F}}}

\newcommand{\QQ}{{\mathbb{Q}}}
\newcommand{\ZZ}{{\mathbb{Z}}}
\newcommand{\NN}{{\mathbb{N}}}

\begin{document}

\title[Rigidity for $F_4(p)$]
{Rational Rigidity for $F_4(p)$}

\date{\today}

\author{Robert Guralnick}
\address{Department of Mathematics, University of
  Southern California, Los Angeles, CA 90089-2532, USA}
  \email{guralnic@usc.edu}

\author{Frank L\"ubeck}
  \address{Lehrstuhl D f\"ur Mathematik, RWTH Aachen,
Pontdriesch 14/16, 52064 Aachen, Germany}
\email{Frank.Luebeck@math.rwth-aachen.de}

\author{Jun Yu}
\address{Beijing International Center for Mathematical Research,
No. 5 Yiheyuan Road, Beijing 100871, China.}
\email{junyu@math.pku.edu.cn}

\thanks{Robert Guralnick was partially supported by the NSF grants
FRG-1265297 and DMS-1302886. He also thanks the Simons Foundation for its
support. This work was initiated in Spring 2013 when Robert Guralnick and
Jun Yu were visiting the Institute for Advanced Study. They thank the
Institute for its hospitality and support.}

\begin{abstract}
We prove the existence of certain rationally rigid triples in
$F_4(p)$ for good primes $p$ (i.e., $p > 3$),
 thereby showing that these groups occur as regular
Galois groups over $\QQ(t)$ and so also over $\QQ$.
 We show that these triples
give rise to rigid triples in the algebraic group and prove that they
generate an interesting subgroup in characteristic $0$.
\end{abstract}

\maketitle

%%%%%%%%%%%%%%%%%%%%%%%%%%%%%%%%%%%%%%%%%%%%%%%%%%%%%%%%%%%%%%%%%%%%%%%%%
\section{Introduction} \label{sec:intro}

The question on which finite groups occur as Galois groups over the field
of rational numbers is still widely open. Even if one restricts to the case
of finite non-abelian simple groups, only rather few types have been realized
as Galois groups over $\QQ$. These include the alternating groups, the
sporadic groups apart from $M_{23}$, and some families of groups of Lie type,
but even over fields of prime order mostly with additional congruence conditions
on the characteristic. In the present paper we show that $F_4(p)$ occur
as Galois groups over $\QQ$ for all good primes $p$.
In \cite{E8}, a similar result was proved by Guralnick and Malle for $E_8(p)$.
A similar result has been proved for $G_2(p)$ for $p > 5$ in \cite{FF} and for 
$p=5$ in \cite{Th}. 

See \cite{E8} and more generally \cite{MM} for more on the general problem
of producing Galois groups over $\QQ$.  In the last few years different
approaches produced Galois groups over $\QQ$ (but not necessarily regular
Galois extensions of $\QQ(t)$; such extensions are
produced by the rigidity method).
See in particular \cite{KLS08, KLS10, Yun, Z13}.

Our proof relies on the well-known rigidity criterion of Belyi, Fried, Matzat 
and Thompson. In addition, we use results about subgroups of a simple linear 
algebraic group of type $\typeF_4$ which contain regular unipotent elements 
(in particular results from \cite{SS97, E8}).

Let $p>3$ be a prime and $k$ be an algebraic closure of the finite field $\FF_p$ 
with $p$ elements. Set $G = F_4(k)$. We consider three conjugacy classes 
$C_1, C_2$ and $C_3$ in $G$. Here $C_1$ is the conjugacy class of involutions 
with centralizer of type $\typeA_1+\typeC_3$, and $C_2$ is the conjugacy class 
of elements of order $2p$ where the $2$-part has a centralizer of type $\typeB_4$ 
and the unipotent part is a long root element in that centralizer. Finally $C_3$ 
is the class of regular unipotent elements.

All classes have representatives in the finite subgroup $F_4(p)$ over the
prime field. And since $p > 3$, all centralizers of elements in these classes
are connected (see section~\ref{Sec-CharComp} for details).
In particular, for each finite subgroup $F_4(q)$, $q$ any power of $p$, the
intersection $C_i(q) := C_i \cap F_4(q)$  ($1 \leq i \leq 3$) is a
single $F_4(q)$-conjugacy class.

Furthermore, the classes $C_i(q)$ are rational, that is for any $g \in C_i(q)$
and $k \in \ZZ$ prime to the order of $g$ we have $g^k \in C_i(q)$. (Let $g
= su = us \in C_i(q)$ with $s$ semisimple and $u$ unipotent. Then $s^k=s$
because the order of $s$ is $1$ or $2$, and $u^k$ is conjugate to $u$ in the
centralizer of $s$ in $F_4(q)$ because the class of $u$ in that centralizer
is uniquely determined by its size.)

We consider the variety \[V(C_1,C_2,C_3) = \{(x,y,z) \in C_1\times C_2
\times C_3 \mid\; xyz = 1\} \subset G^3,\] on which $G$ acts by componentwise 
conjugation.

Our main result is as follows.

\begin{thm} \label{thm:rigidity}
 Let $p$, $k$, $G=F_4(k)$, $X := V(C_1,C_2,C_3)$ be as above.
 Then $G$
 has a single regular orbit on $X$ and if $(x, y, z) \in X$,
 then $\langle x, y \rangle \cong F_4(p)$.   Moreover, for any power $q$ of
 $p$ the subgroup $F_4(q)$ has a single regular orbit on $X(q)$, the
 $\FF_q$-rational points of $X$.
\end{thm}

As noted the rigidity criterion as stated in~\cite[I.~Thm.~4.8]{MM} implies 
the following (we use that the center of $F_4(p)$ is trivial, and that the 
classes $C_i(p)$ are rational). 

\begin{cor}  \label{cor:Galois}  If $p >3 $ is a prime, then $F_4(p)$ is a regular 
Galois group over $\QQ(t)$ and in particular is a Galois group over $\QQ$.
\end{cor}

\textbf{Remark.} This result was already obtained for $p\geq 19$ and
$p\equiv 2, 6, 7 \textrm{ or }11 \pmod{13}$ in~\cite{MaExc}, using a
similar argument for a different triple of conjugacy classes.

\smallskip

As usual for applications of the rigidity criterion our proof has two main
parts. In section~\ref{Sec-Gen} we show that any triple in $V(C_1,C_2,C_3)$
generates a subgroup which contains a conjugate of $F_4(p)$. We use
non-trivial results on maximal subgroups of $G$ containing a regular
unipotent element and some representation theory for Levi subgroups of $G$.

In section~\ref{Sec-CharComp} we compute a certain estimate for the structure
constants of the class triples $C_i(q)$ ($1 \leq i \leq 3$) for all powers
$q$ of $p$.
We use Deligne-Lusztig theory for these computations. These computations are
a bit more complicated than, e.g., in~\cite{E8} because many character
values on the mixed class $C_2(q)$ are needed.

In section~\ref{Sec-Proof} we combine the results and prove
Theorem~\ref{thm:rigidity}.

Since $G=F_4(k)$ has a single regular orbit
on $V(C_1,C_2,C_3)$ for $k$ algebraically
closed of  good positive characteristic, it follows that the same is true if
$k$ is an algebraically closed field of
characteristic $0$.  Precisely as in \cite{E8} since $C_G(z)$ for $z \in
C_3$ is an affine
space of dimension $4$, it follows that this torsor is trivial and so
$G$ has a single orbit on the $k$-points of the variety described above
over any field of characteristic $p \ne 2, 3$ (including characteristic zero).

Arguing exactly as in \cite{E8},
we can also produce
such triples over $G(\ZZ_p)$ and so deduce:

\begin{thm} \label{thm:char 0}
 Let $k$ be an algebraically closed field of characteristic $0$.
 Let $G=F_4(k)$ and $C_i$, $1\leq i\leq 3$, be the conjugacy classes
 described above. Let $X := V(C_1,C_2,C_3)$.
 For a triple $c \in X$,
 let $\Gamma(c)$ denote the group generated by $c$.
 \begin{enumerate}
  \item  For any $c \in X$, $\Gamma(c)$ is Zariski dense in $G$.
  \item  If $k_0$ is a subfield of $k$, then $X(k_0)$, the $k_0$-rational
  points of $X$, is a single
   $F_4(k_0)$-orbit (where $F_4(k_0)$ is the split group over $k_0$).
  \item Set $R=\ZZ[1/6]$. There exists $c\in X(R)$ such that $\Gamma(c)\le F_4(R)$ 
  and surjects onto $F_4(R/pR)$ for any good prime $p$. In particular, $\Gamma(c)$ 
  is dense in $F_4(\ZZ_p)$ with respect to the $p$-adic topology for any good 
  prime $p$.
 \end{enumerate}
\end{thm}

We thank Zhiwei Yun for suggesting the conjugacy
classes considered in this paper (see also the discussion in~\cite[7.]{E8}).

\section{A Generation Result}\label{Sec-Gen}

Fix a prime $p > 3$ and let $k$ be an algebraic closure of $\FF_p$.
Let $G=F_4(k)$ and let $C_i$, $1 \leq i \leq 3$, be the conjugacy classes
of $G$ described in the introduction.

\begin{thm}\label{T:generation}
If $(x,y,z)\in V(C_1,C_2,C_3)$,
then $H:=\langle x,y,z\rangle$ contains a conjugate of
$F_4(p) \leq F_4(k)$.
\end{thm}

\noindent\textbf{Remark.} In section~\ref{Sec-Proof} we will show, using
theorem~\ref{thmstructconst}, that in fact $H$ is a conjugate of $F_4(p)$.

\smallskip

\begin{proof} (of the theorem)

(1) Let $y=y_sy_u=y_uy_s$ with $y_s$ semisimple and $y_u$ unipotent be the
Jordan decomposition of $y$. Since $y_s$ and $y_u$ are powers of $y$ we know
that $H$ contains two involutions $x$ and $y_s$ which are not conjugate and
two unipotent elements $y_u$ and $z$ which are not conjugate in $ F_4(k)$.

(2) The finite maximal subgroups of $F_4(k)$ which contain a regular
unipotent element which do not contain a conjugate of $F_4(p)$ and are also
not contained in any proper positive dimensional subgroup
(called \emph{Lie primitive} subgroups) were determined in~\cite[Theorem
3.4(a)]{E8}. These subgroups have the property that $p'$-powers of any
non-trivial unipotent element yield representatives of all their non-trivial
unipotent conjugacy classes. Hence all unipotent
elements in these subgroups are regular unipotent in $F_4(k)$ and so they
cannot contain $H$.

(3) The maximal closed subgroups of $F_4(k)$ of positive dimension
which contain a regular unipotent element where determined
in~\cite[Theorem A]{SS97}. These are parabolic or reductive of type $\typeA_1$
(note that we assume $p>3$). It is clear that $H$ cannot be contained in a
subgroup of type $\typeA_1$ which contains only one class of unipotent elements.
In the remainder of this proof we show that $H$  also cannot be contained in
any proper parabolic subgroup of $G$.

(4) Let $P < G$ be a maximal parabolic subgroup and assume $H =\langle
x,y,z\rangle \leq P$. Let $U \lhd P$ be the unipotent radical of $P$
and $L$ a Levi complement. So, $P = U \rtimes L$.
Every element $g \in P$ can be uniquely written as $g = g' \bar{g}$ with $g'
\in U$ and $\bar{g} \in L$ and we have a homomorphism $P \to L$, $g \mapsto
\bar{g}$, which maps semisimple elements to semisimple elements and unipotent
elements to unipotent elements.

For the involution $x = x' \bar{x}$ we get that $\bar{x}$ also has order $2$
and that $x \bar{x} = x'$ is unipotent and so of odd order. Hence $x$ and
$\bar{x}$ are conjugate in the dihedral group they generate and so in $G$.
The same argument shows that $y_s$ and $\bar{y}_s$ are conjugate involutions
in $G$.

To find the possibilities for $\bar{y}_u$ we use the fact that for any
unipotent element $u \in P$ the element $\bar{u}$ is contained in the
closure of the $G$-conjugacy class $u^G$, see~\cite[(9.5.2)]{Lu84}.
Using this and~\cite[13.4]{Ca}, it follows that $\bar{y}_u$ is trivial
or is conjugate to a long root element in $G$.

The relation $xyz=1$ implies $\bar{x}\bar{y}_s\bar{y}_u\bar{z} = 1$.

Our strategy is to find a representation of $L$ such that we can show
that $\bar{x}\bar{y}_s\bar{y}_u$ has an eigenvalue $-1$ which is in
contradiction to $\bar{z}^{-1}$ being unipotent. Note that in fact
$\bar{z}$ is a regular unipotent element in $L$ (we will not use this
though).

(5) We use the setup described in~\cite[Section 2]{BL13} to determine
classes of involutions in $G$ and some Levi subgroups, and to compute the
eigenvalues of involutions in various representations.

The group $G$ can be described by $k$ and a root datum
$(X,R,Y,R^\vee)$. Here $X$ and
$Y$ are dual $\ZZ$-lattices of rank 4. As $\ZZ$-basis of $Y$ we choose a set
of simple coroots in $R^\vee$ (note that $G$ is simply-connected).
The dual basis of $X$ consists of the
fundamental weights of $G$ and the corresponding set of simple roots with
respect to this basis is given by the transposed Cartan matrix of the root
system of $G$, that is the rows of
{\footnotesize $\left(\begin{array}{rrrr}%
2&-1&0&0\\%
-1&2&-2&0\\%
0&-1&2&-1\\%
0&0&-1&2\\%
\end{array}\right)%
$} (the first two simple roots are long).
From this it is straightforward to write down matrices for the Coxeter
generators of the Weyl group $W$ of $G$ acting on $X$:
$s_\alpha(x) = x - \langle x,
\alpha^\vee \rangle \alpha$, where $\alpha$ is a simple root, $\alpha^\vee$
the corresponding coroot, and $\langle.,.\rangle$ is the pairing between $X$
and $Y$. Similarly, matrices for the Weyl group generators acting on $Y$ can
be computed.

This also describes the action of the Weyl group $W$ on a maximal torus of
$G$ which is isomorphic to $T = Y \otimes_\ZZ (\QQ/\ZZ)_{p'}$.
A torus element $t = (t_1,t_2,t_3,t_4) \in T$ can be evaluated at a weight
$x = (x_1,x_2,x_3,x_4) \in X$ by $\sum_{i=1}^4 x_it_i \in (\QQ/\ZZ)_{p'}
\cong k^\times$. The roots which have $t$ in their kernel form the root
system of the centralizer of $t$.

Conjugacy classes of semisimple elements are in bijection with $W$-orbits in
$T$ (see~\cite[3.7]{Ca}). Computing the $W$-orbits on the 15 elements of
order 2 in $T$, we recover the 2 conjugacy classes of involutions in
$G$ (the class of $x$ with centralizer of type $\typeA_1 + \typeC_3$
and the class of $y_s$ with centralizer of type $\typeB_4$) and all their
representatives in $T$. We will use this below to find the fusion of classes
of involutions in Levi subgroups into $G$.

Let $L_i$ ($1 \leq i \leq 4$) be the standard Levi subgroups of $G$, whose
root datum we get from the root datum of $G$ by removing the $i$-th simple
root and the corresponding coroot. We denote by $P_i$ the corresponding maximal
standard parabolic subgroup and $U_i$ its unipotent radical.
The derived subgroup $L_i'$ of $L_i$ is also of simply-connected type and
the intersection of $T \cap L_i'$ consists of the $t = (t_1,t_2,t_3,t_4)$
with $t_i=0$ (see~\cite[2.6]{BL13}). We now consider each $P_i$
separately\footnote{Using an approach in~\cite[Section 3.1]{BFM}, we can
describe the isomorphism type of each $L_{i}$ ($1\leq i\leq 4$). Precisely
we have \[L_1\cong(\Sp_{6}(k)\times k^{\ast})/\langle(-I,-1)\rangle,\]
\[L_2\cong L_3\cong(\SL_{3}(k)\times\SL_{2}(k)\times
k^{\ast})/\langle(\omega I, I,\omega^{-1}), (I,-I,-1)\rangle,\] \[L_4\cong
(\Spin_{7}(k)\times k^{\ast})/\langle(-1,-1)\rangle.\] By \cite[Table
2]{HY13}, the centralizer of $x$ is a subgroup isomorphic to
$(\Sp_{6}(k)\times\SL_{2}(k))/\langle(-I,-I)\rangle$, which contains a
conjugate of $L_1$; the centralizer of $y_{s}$ is a subgroup
isomorphic to $\Spin_{9}(k)$, which contains a conjugate of $L_4$. Then,
using the description of fusions of involutions in these symmetric
subgroups (pages 414-415 in \cite[Table 2]{HY13}), we can also identify
the possible classes of $\bar{x}$ and $\bar{y}_{s}$ in $L_1$ (or $L_4$).
By describing the fusion of involutions in a connected semisimple subgroup
of $G$ which is isomorphic to 
$(\SL_{3}(k)\times\SL_{3}(k))/\langle(-I,-I)\rangle$, we can also identify 
the possible classes of $\bar{x}$ and $\bar{y}_{s}$ in $L_2$ (or $L_3$).}

(6) Case $P_1$. The commutator subgroup $L_1'$ is isomorphic to $\Sp_6(k)$.
We consider $\bar{x}\bar{y}_s\bar{y}_u\bar{z} = 1$. The only
unipotent class of $L_i$ (or $L_i'$) which is conjugate in $G$ to a long
root element is the class of a long root element in $L_1'$ (see the
description of unipotent classes in~\cite[13.1]{Ca} which also holds in
characteristic $>3$). So, $\bar{y}_u$ is trivial or a long root element in
$L_1'$.

Using the Weyl group of $L_1$ to compute its orbits on the involutions in
$T$ we see that $L_1$ has one class of involutions which are conjugate to
$y_s$ in $G$ and three classes of involutions which are conjugate to $x$ in
$G$. Only one of the classes conjugate to $x$ is not contained in $L_1'$.
Since $\bar{x}\bar{y}_s\bar{y}_u\bar{z} = 1$ and $\bar{y}_{s}$ and the unipotent
elements $\bar{y}_u$ and $\bar{z}$ are in $L_1'$ it follows that $\bar{x}$
must also lie in $L_1'$.

We consider the irreducible representation of $L_1$ with highest weight
$(0,0,0,1)$. This is an extension of the 6-dimensional natural
representation of $L_1'$ to $L_1$ (compare~\cite[3.5 and 3.6(a)]{BL13}).
It weights consist of the $W$-orbit of $(0,0,0,1)$.

The image of $\bar{y}_u$ is trivial or has a unique non-trivial Jordan block
of size $2$; in any case $\bar{y}_u$ has at least a 5-dimensional
$1$-eigenspace. From the weights and class representatives we compute that
$\bar{y}_s$ has a 2-dimensional $1$-eigenspace and a 4-dimensional
$(-1)$-eigenspace. For $\bar{x}$ there are two possibilities, it has
eigenvalue $(-1)$ either on the full space or on a 2-dimensional subspace.
In any case we find that $\bar{x}\bar{y}_s$ has at least a 2-dimensional
$(-1)$-eigenspace and so $\bar{x}\bar{y}_s\bar{y}_u$ has at least a
1-dimensional $(-1)$-eigenspace. On the other hand $\bar{z}$ is unipotent
and has no eigenvalue $(-1)$. This shows that $H$ cannot be contained in
$P_1$.

(7) Case $P_2$. Here $L_2'$ is isomorphic to $\SL_2(k) \times \SL_3(k)$ and
only the nontrivial unipotent elements in the $\SL_2(k)$-component are
conjugate to $y_u$ in $G$. We consider the representation of $L_2$ with highest
weight $(1,0,0,0)$ which is an extension of the natural
2-dimensional representation of the $\SL_2(k)$-component of $L_2'$.

The group $L_2$ contains only one class of involutions conjugate to $y_s$ in $G$, 
this is contained in $L_2'$ and its elements act as identity on the 2-dimensional
representation. Since $\bar{x}\bar{y}_s\bar{y}_u\bar{z} = 1$ and $\bar{y}_{s}$ 
and the unipotent elements $\bar{y}_u$ and $\bar{z}$ are in $L_2'$ it follows 
that $\bar{x}$ must also lie in $L_2'$.

There are four classes of $L_2$ which belong to the class of
$x$ in $G$, two of them are in $L_2'$ and so may contain $\bar{x}$.
In both cases $\bar{x}$ acts as $(-1)$ times the identity on the
2-dimensional representation. We get that the image of $\bar{x}\bar{y}_s
\bar{y}_u$ has eigenvalue $(-1)$ in this 2-dimensional representation.
This is not possible because $\bar{z}$ is unipotent.

(8) Case $P_3$. Here $L_3'$ is isomorphic to $\SL_3(k) \times \SL_2(k)$
and only the root elements of the $\SL_3(k)$-component are conjugate to
$y_u$ in $G$. We consider the representation of $L_3$ with highest weight
$(1,0,0,0)$  which
extends the 3-dimensional natural representation of the $\SL_3(k)$-component
of $L_3'$. We find that for all possibilities $\bar{y}_s$ acts trivially in
this representation and the image of $\bar{x}$ has a
$(-1)$-eigenspace of dimension $2$. Furthermore, the image of $\bar{y}_u$
has at most one Jordan block of size $2$. This leads to a contradiction as
in the previous cases.

(9) Case $P_4$. Now $L_4'$ is isomorphic to a spin group $\Spin_7(k)$. We
proceed as in the other cases, now using the 7-dimensional representation
with highest weight $(1,0,0,0)$ which is an extension of the 7-dimensional
representation of $L_4'$ with image $\SO_7(k)$ (the weights are $0$ and
the orbit of the highest weight). The image of a long root element in
$L_4'$ has two non-trivial Jordan blocks of size $2$ in this representation.
The rest of the argument is similar to the other cases.

We have shown that $H$ is not contained in any (maximal) parabolic subgroup.
\end{proof}

\section{Character Computations}\label{Sec-CharComp}

As in the introduction, let $p>3$ be a prime, $k$ be an algebraic closure
of $\FF_p$ and $G = F_4(k)$ be a simple linear algebraic
group with root system of type $\typeF_4$. Recall that we have fixed three
conjugacy classes $C_1$, $C_2$ and $C_3$ of $G$ as follows. The class $C_1$ 
consists of involutions whose centralizer is a connected reductive subgroup 
with root system of type $\typeA_1 + \typeC_3$.
Note that $G$ is of simply-connected type and so
all centralizers of semisimple elements are connected (\cite[3.5.6]{Ca}).
The class $C_2$ contains elements with Jordan decomposition $su=us$ where
$s$ is an involution with a reductive centralizer with root system of type
$\typeB_4$ and $u$ is a long root element in the centralizer $C_G(s)$ of $s$.
More precisely, $C_G(s)$ is isomorphic to the group $\Spin_9(k)$ and
the image of $u$ in the representations as $\SO_9(k)$ has Jordan block sizes
$2^21^5$. The centralizer of $u$ in $C_G(s)$, which is the centralizer
$C_G(su)$, is connected (\cite[14.3]{Lu84}). The third class $C_3$ consists
of the regular unipotent elements of $G$. Since we are in good
characteristic for $G$ the centralizer in $G$ of $u \in C_3$ is also
connected (\cite[5.6.1]{Ca}).

Now we consider for each power $q$ of $p$ the finite group $F_4(q)$. Since
our conjugacy classes $C_1$, $C_2$ and $C_3$ contain elements in $F_4(p)$
and since these elements have connected centralizer in $G$ it follows from
the Lang-Steinberg theorem that $C_1$, $C_2$ and $C_3$ intersect each $F_4(q)$
in a single $F_4(q)$-conjugacy class.

Despite a lot of theory that is known about the ordinary characters of
the groups $F_4(q)$ it is still impossible to compute their full character
tables. However, using partial information about the character table
we can work out the following result.

\begin{thm}\label{thmstructconst}
Let $X = V(C_1,C_2,C_3)$ be as before and $X(q)$ be the set of
$\FF_q$-rational points of $X$ for any power $q$ of $p$. There exists a
positive constant $c$ such that for all $p$ and $q$:
\[ 1-\frac{c}{q} < \frac{|X(q)|}{|F_4(q)|} < 1+\frac{c}{q}.\]
\end{thm}

\noindent\textbf{Remark.}
In particular $X(q)$ is not empty for large enough $q$. In
section~\ref{Sec-Proof} we will show that
in fact the quotient $|X(q)| / |F_4(q)|$ is precisely $1$ for all $q$.

\smallskip

\begin{proof} (of the theorem)

(1) Let $x \in C_1 \cap F_4(p)$, $y \in C_2 \cap F_4(p)$ and
$z \in C_3 \cap F_4(p)$.
The quantity $|X(q)| / |F_4(q)|$ can be computed as class structure
constant by the following formula (see for example~\cite[I. Thm. 5.8]{MM}):
\[ \frac{|X(q)|}{|F_4(q)|} = \frac{|F_4(q)|}{c_x(q) c_y(q) c_z(q)}
\sum_{\chi \in {\rm Irr}(F_4(q))} \frac{\chi(x) \chi(y) \chi(z)}{\chi(1)},\]
where $c_x(q)$, $c_y(q)$ and $c_z(q)$ are the centralizer orders of
$x, y, z$, respectively,  in $F_4(q)$.

The factor $|F_4(q)|/(c_x(q) c_y(q) c_z(q))$ is a quotient of two monic
polynomials in  $q$ of the same degree, so it tends to $1$ for growing $q$.
It remains to show a lower and an upper
bound as in the theorem for the sum in the formula.

(2) Since $G$ has trivial (and hence connected) center we can use Lusztig's
parameterization of the irreducible characters of $F_4(q)$
in~\cite[Ch.4]{LuBook}. Here, to each conjugacy class of semisimple elements
in the dual group of $F_4(q)$, which is isomorphic to $F_4(q)$ itself,
Lusztig associates a subset of the irreducible characters, which we now call
a \emph{Lusztig series}.  For a semisimple element $s$ of the dual group
the corresponding Lusztig series is parameterized by some set which only
depends on the twisted Lie type of the centralizer of $s$. We do not need
the full details of this description because of the next two remarks.

(3) We only need to know character values of irreducible characters which
are non-zero on all three classes $C_1$, $C_2$ and $C_3$. The characters
which have non-zero value on class $C_3$ (regular unipotent elements) are
called \emph{semisimple characters}. Each Lusztig series contains exactly
one semisimple character, see~\cite[8.3, 8.4]{Ca} (note that we use here
our assumption that $p$ is a good prime for $G$).

(4) Semisimple characters are uniform, which means that they are linear
combinations of Deligne-Lusztig characters. An explicit description of
this linear combination is given in~\cite[8.4.6]{Ca}.

(5) The semisimple conjugacy classes of $F_4(q)$ (and its isomorphic dual
group) were determined in~\cite{Sh74}. But we need a slightly more precise
description.
We first describe representatives of the $F_4(q)$-classes of centralizers
of semisimple elements in the dual group (this does not depend on $q$) 
and then for each of them the elements in the center of this centralizer. 
We also (re)-compute the number of semisimple classes whose elements have 
centralizer in each class.

(6) We need representatives of the $F_4(q)$-classes of maximal tori
in $G$ and parameterizations of their elements, this can be done as in~(5).
Furthermore, if $T$ and $T^*$ are dual maximal tori in $G$ and its dual
group, stable under a Frobenius morphism over $\FF_q$, we need to associate
to a pair $(T^*,s)$ with $s \in T^*(q)$ an explicit linear character of
$T(q)$. This is described in the proof of~\cite[13.7]{DM91}. Furthermore,
we identify representatives of the semisimple parts of elements in our
classes $C_1$ and $C_2$ in these tori.

(7) We compute the values of Deligne-Lusztig characters on the trivial class
and on our classes $C_1$, $C_2$ and $C_3$ via the character
formula~\cite[7.2.8]{Ca}. In addition to~(6) we need some values of Green
functions for this computation. For the trivial class and $C_1$ these can be
computed from certain torus orders, see~\cite[7.5.1]{Ca}. On the regular
unipotent elements Green functions have value $1$, see~\cite[8.4.3]{Ca}.
And on the long root elements in groups of type $\typeB_4$ they can be computed
using~\cite{Sh87}.

(8) Now we can compute the values of all semisimple characters on the
classes we consider as linear combinations of certain Deligne-Lusztig
characters using~(7),~(6) and~(4). We get parameterized descriptions of the
values for all Lusztig series corresponding to semisimple elements $s$ with
a fixed class of centralizer. We call such a subset of semisimple characters
a \emph{semisimple character type}. These values are polynomials in $q$ whose
coefficients are sums of powers of $-1$ (the root of unity of order $2$)
where these powers depend on certain parameters which distinguish the
characters within a character type.

For more details on such computations we refer to~\cite[Section 2]{BL13}
and~\cite{FLDiss}. We have used computer programs which are based on
code provided by the CHEVIE package~\cite{CH}. The results are a bit
different according to the cases $q \equiv 1, 5, 7, 11 \pmod{12}$,
so all computations must be done separately for each of these congruence
classes of $q$ modulo $12$.

(9) In the statement of the theorem the $1$ on the right hand side
comes from the trivial character in the sum.  To show that all other
characters contribute something of absolute value less than a (big enough)
constant times $1/q$, we can show this for each semisimple
character type corresponding to non-trivial semisimple elements in the
dual group separately (there are about $100$ semisimple character types).

Let $\chi_s$ be the semisimple character corresponding to a semisimple
element $s$ in the dual group, and let $n_s$ be the number of characters
in the semisimple character type of $s$ ($n_s$ is also a polynomial
evaluated at $q$).
For most semisimple character types we find that the degrees as polynomials
in $q$ of $f_s := n_s \chi_s(x) \chi_s(y) \chi_s(z)
|F_4(q)| / \chi_s(1)$ and
of $|F_4(q)|$ are the same. So, just comparing the degrees of these
polynomials is not good enough, we must show that at least the leading
coefficients of the $f_s$ sum up to zero when we sum over all $s$
corresponding to a fixed semisimple character type.

We illustrate this with an example. Consider elements $s$ which are regular
elements in a maximal torus of order $(q+1)(q^3+q^2+q+1)$. These are
parameterized by two integers $1 \leq m \leq q^3+q^2+q+1$ and $1 \leq n
\leq q+1$ with the exception of $O(q^3)$ parameter pairs (which correspond
to non-regular elements in that torus).
The leading coefficient of $\chi_s(x) \chi_s(y) \chi_s(z)
|F_4(q)| / \chi_s(1)$ is $((-1)^{n+m} + (-1)^n) q^{48}$. Since $q$ is odd we
see that the coefficients of $q^{48}$ cancel out if we sum over all
parameters $m, n$ as above. Hence the complete summation over the characters
in this semisimple character type can be bounded by a constant times
$q^{48+3}$ while $|F_4(q)| = q^{52}+O(q^{51})$.
\end{proof}

\section{Proof of Theorem 1.1}\label{Sec-Proof}

Let $k$ be an algebraic closure of $\FF_p$, $p > 3$, and $G = F_4(k)$.
Let $X = \{(x,y,z)\mid\; x \in C_1, y \in C_2, z \in C_3, xyz = 1\}$ be as
before.

In Theorem~2.1 we have seen that the group generated by any triple in $X$
contains a conjugate of $F_4(p)$. Therefore, the stabilizer of each triple
must centralize a conjugate of $F_4(p)$, hence it is trivial (and connected).
This shows that each $G$-orbit on $X$ is regular.

Each triple in $X$ is contained in a finite subgroup $F_4(q)$ for some power
$q$ of $p$ (because $k = \bar \FF_p$), and so its $G$-orbit is stable under
all Frobenius morphisms of $G$ over $\FF_{q^m}$, $m \in \NN$. It follows from
the Lang-Steinberg theorem that the $\FF_{q^m}$-rational points of such an
orbit form one regular $F_4(q^m)$-orbit.

In Theorem~3.1 we have computed an estimate for $|X(q)|/|F_4(q)|$ which
shows that for all sufficiently large $q$ the set $X(q)$ is not empty
and consists of exactly one $F_4(q)$-orbit. Hence $X$ consists of a
single $G$-orbit.

Now $X$ is invariant under the Frobenius morphism of $G$ over the prime
field $\FF_p$ (because our classes $C_i$ are defined over $\FF_p$)) and 
so over any $\FF_q$ for powers $q$ of $p$. And since it is a single $G$-orbit
we can apply again the Lang-Steinberg theorem and see that the set
$X(q)$ of $\FF_q$-rational points of $X$ is not empty and is a single
$F_4(q)$-orbit. Using this for $q=p$ and Theorem~2.1 again a triple from
$X(p)$ must generate $F_4(p)$.
Therefore, all triples from $X$ generate conjugates of $F_4(p)$.


\begin{thebibliography}{999999}

\bibitem {BFM} A.~Borel, R.~Friedman, J.~Morgan, \emph{Almost
commuting elements in compact Lie groups.} Mem. Amer. Math. Soc. 157
(2002), no. 747.

\bibitem {BL13}{\sc O. Brunat, F. L\"ubeck}, On defining
characteristic representations of finite reductive groups. J. Algebra 
{\bf 395} (2013), 121--141.

\bibitem {Ca}{\sc R. W. Carter}, \emph{Finite Groups of Lie Type.
Conjugacy Classes and Complex Characters}. Wiley Classics Library. John
Wiley \& Sons, Chichester, 1993.

\bibitem{DM91}{\sc F. Digne, J. Michel}, \emph{Representations
of Finite Groups of Lie Type}. LMS Student Texts, {\bf 21}. Cambridge
University Press, Cambridge, 1991.

\bibitem{FF} {\sc W. Feit, P. Fong} Rational rigidity of $G_2(p)$
for any prime $ p>5$, Proceedings of the Rutgers group theory year,
1983--1984 (New Brunswick, N.J., 1983--1984), 323--326, Cambridge Univ.
Press, Cambridge, 1985.

\bibitem{CH}{\sc M.~Geck, G.~Hiss, F.~L{\"u}beck, G.~Malle, and
G.~Pfeiffer}, C{HEVIE}---a system for computing and processing generic
character tables. Appl. Algebra Engrg. Comm. Comput. {\bf 7(3)} (1996),
175--210.

\bibitem{E8}{\sc R. M. Guralnick, G. Malle}, Rational rigidity for
$E_8(p)$. Compos. Math. {\bf 150} (2014), 1679--1702.

\bibitem{HY13}{\sc J.-S.~Huang, J.~Yu}, Klein four-subgroups of Lie
algebra automorphisms. Pacific J. Math. {\bf 262} (2013), no. 2, 397--420.

\bibitem{KLS08}{\sc C. Khare, M. Larsen, G. Savin}, Functoriality
and the inverse Galois problem. Compos. Math. {\bf 144} (2008), 541--564.

\bibitem{KLS10}{\sc C. Khare, M. Larsen, G. Savin}, Functoriality
and the inverse Galois problem. II. Groups of type $\typeB_n$ and
$\typeG_2$. Ann.  Fac. Sci. Toulouse Math. (6) {\bf 19} (2010), 37--70.

\bibitem{FLDiss}{\sc F. L\"ubeck}, Charaktertafeln f\"ur die
Gruppen $\CSp_{6}(q)$ mit ungeradem $q$ und $\Sp_6(q)$ mit geradem $q$.
Dissertation, Universit\"at Heidelberg (1993).

\bibitem{Lu84}{\sc G. Lusztig}, Intersection cohomology on a
reductive group. Invent. Math. {\bf 75} (1984), 205--272.

\bibitem{LuBook}{\sc G. Lusztig}, \emph{Characters of reductive
groups over a finite field.} Annals of Mathematics Studies {\bf 107},
Princeton University Press, 1984.

\bibitem{MaExc}{\sc G. Malle}, Exceptional groups of {L}ie type as
Galois groups. J. Reine Angew. Math. {\bf 392} (1988), 70--109.

\bibitem{MM}{\sc G. Malle, B. H. Matzat}, \emph{Inverse Galois
Theory.} Springer Monographs in Mathematics. Springer-Verlag, Berlin,
1999.

\bibitem{SS97}{\sc J. Saxl, G. M. Seitz}, Subgroups of algebraic
groups containing regular unipotent elements. J. London Math. Soc. {\bf55}
(1997), 370--386.

\bibitem{Sh74}{\sc T. Shoji}, The conjugacy classes of Chevalley
groups of type ($\typeF_4$) over finite fields of characteristic $p \neq 2$. J.
Fac. Sci. Univ. Tokyo {\bf 21} (1974), 1--17.

\bibitem{Sh87}{\sc T. Shoji}, Green functions of reductive groups
over a finite field. Proc. Symp. Pure Math. {\bf 47} (1987), 289--301.

\bibitem{Th} {\sc J. Thompson} Rational rigidity of $G_2(5)$,
Proceedings of the Rutgers group theory year, 1983--1984 (New Brunswick,
N.J., 198--1984), 321--322, Cambridge Univ. Press, Cambridge, 1985.

\bibitem{Yun}{\sc Z. Yun}, Motives with exceptional Galois
groups and the inverse Galois problem. Invent. Math. {\bf 196} (2014), no.
2, 267--337.

\bibitem{Z13}{\sc D. Zywina}, The inverse Galois problem for
$\mathrm{\PSL_2}(\FF_p)$, Preprint, arXiv:1303.3646v1.


%\bibitem {BFM} A.~Borel, R.~Friedman, J.~Morgan, \emph{Almost
%commuting elements in compact Lie groups.} Mem. Amer. Math. Soc. 157
%(2002), no. 747.

%\bibitem {BL13}{\sc O. Brunat, F. L\"ubeck}, On defining
%characteristic representations of finite reductive groups. J. Algebra
%{\bf 395} (2013), 121--141.

%\bibitem {Ca}{\sc R. W. Carter}, \emph{Finite Groups of Lie Type.
%Conjugacy Classes and Complex Characters}. Wiley Classics Library. John
%Wiley \& Sons, Chichester, 1993.

%\bibitem[DM91]{DM91}{\sc F. Digne, J. Michel}, \emph{Representations
%of Finite Groups of Lie Type}. LMS Student Texts, {\bf 21}. Cambridge
%University Press, Cambridge, 1991.

%\bibitem[FF85]{FF} {\sc W. Feit, P. Fong} Rational rigidity of $G_2(p)$
%for any prime $ p>5$, Proceedings of the Rutgers group theory year,
%1983--1984 (New Brunswick, N.J., 1983--1984), 323--326, Cambridge Univ.
%Press, Cambridge, 1985.

%\bibitem[CH94]{CH}{\sc M.~Geck, G.~Hiss, F.~L{\"u}beck, G.~Malle, and
%G.~Pfeiffer}, C{HEVIE}---a system for computing and processing generic
%character tables. Appl. Algebra Engrg. Comm. Comput. {\bf 7(3)} (1996),
%175--210.

%\bibitem[GM12]{E8}{\sc R. M. Guralnick, G. Malle}, Rational rigidity for
%$E_8(p)$. Compos. Math. {\bf 150} (2014), 1679--1702.

%\bibitem[HY13]{HY13}{\sc J.-S.~Huang, J.~Yu}, Klein four-subgroups of Lie
%algebra automorphisms. Pacific J. Math. {\bf 262} (2013), no. 2, 397--420.

%\bibitem[KLS08]{KLS08}{\sc C. Khare, M. Larsen, G. Savin}, Functoriality
%and the inverse Galois problem. Compos. Math. {\bf 144} (2008), 541--564.

%\bibitem[KLS10]{KLS10}{\sc C. Khare, M. Larsen, G. Savin}, Functoriality
%and the inverse Galois problem. II. Groups of type $\typeB_n$ and
%$\typeG_2$. Ann.  Fac. Sci. Toulouse Math. (6) {\bf 19} (2010), 37--70.

%\bibitem[FL93]{FLDiss}{\sc F. L\"ubeck}, Charaktertafeln f\"ur die
%Gruppen $\CSp_{6}(q)$ mit ungeradem $q$ und $\Sp_6(q)$ mit geradem $q$.
%Dissertation, Universit\"at Heidelberg (1993).

%\bibitem[Lu84]{Lu84}{\sc G. Lusztig}, Intersection cohomology on a
%reductive group. Invent. Math. {\bf 75} (1984), 205--272.

%\bibitem[Lu84b]{LuBook}{\sc G. Lusztig}, \emph{Characters of reductive
%groups over a finite field.} Annals of Mathematics Studies {\bf 107},
%Princeton University Press, 1984.

%\bibitem[M88]{MaExc}{\sc G. Malle}, Exceptional groups of {L}ie type as
%Galois groups. J. Reine Angew. Math. {\bf 392} (1988), 70--109.

%\bibitem[MM99]{MM}{\sc G. Malle, B. H. Matzat}, \emph{Inverse Galois
%Theory.} Springer Monographs in Mathematics. Springer-Verlag, Berlin,
%1999.

%\bibitem[SS97]{SS97}{\sc J. Saxl, G. M. Seitz}, Subgroups of algebraic
%groups containing regular unipotent elements. J. London Math. Soc. {\bf55}
%(1997), 370--386.

%\bibitem[Sh74]{Sh74}{\sc T. Shoji}, The conjugacy classes of Chevalley
%groups of type ($\typeF_4$) over finite fields of characteristic $p \neq 2$. J.
%Fac. Sci. Univ. Tokyo {\bf 21} (1974), 1--17.

%\bibitem[Sh87]{Sh87}{\sc T. Shoji}, Green functions of reductive groups
%over a finite field. Proc. Symp. Pure Math. {\bf 47} (1987), 289--301.

%\bibitem[Th85]{Th} {\sc J. Thompson} Rational rigidity of $G_2(5)$,
%Proceedings of the Rutgers group theory year, 1983--1984 (New Brunswick,
%N.J., 198--1984), 321--322, Cambridge Univ. Press, Cambridge, 1985.

%\bibitem[Yu12]{Yun}{\sc Z. Yun}, Motives with exceptional Galois
%groups and the inverse Galois problem. Invent. Math. {\bf 196} (2014), no.
%2, 267--337.

%\bibitem[Z13]{Z13}{\sc D. Zywina}, The inverse Galois problem for
%$\mathrm{\PSL_2}(\FF_p)$, Preprint, arXiv:1303.3646v1.


\end{thebibliography}
\end{document}